\documentclass[11pt]{article}
\usepackage{amsmath,amssymb,amsthm,bm,graphicx,color,epsfig}

\newtheorem{theorem}{Theorem}[section]

\newtheorem{definition}[theorem]{Definition}

\newcommand{\R}{\mathbb{R}}
\renewcommand{\S}{\mathbb{S}}
\renewcommand{\i}{\imath}

\newcommand{\eps}{\varepsilon}
\newcommand{\p}{\partial}

\begin{document}

\title{Fast Directional Computation for the High Frequency Helmholtz
  Kernel in Two Dimensions}
\author{Bj\"{o}rn Engquist and Lexing Ying\\
  Department of Mathematics, University of Texas, Austin, TX 78712}

\date{February 2008}
\maketitle

\begin{abstract}
  This paper introduces a directional multiscale algorithm for the two
  dimensional $N$-body problem of the Helmholtz kernel with
  applications to high frequency scattering. The algorithm follows the
  approach in [Engquist and Ying, SIAM Journal on Scientific
  Computing, 29 (4), 2007] where the three dimensional case was
  studied. The main observation is that, for two regions that follow a
  directional parabolic geometric configuration, the interaction
  between the points in these two regions through the Helmholtz kernel
  is approximately low rank. We propose an improved randomized
  procedure for generating the low rank representations.  Based on
  these representations, we organize the computation of the far field
  interaction in a multidirectional and multiscale way to achieve
  maximum efficiency. The proposed algorithm is accurate and has the
  optimal $O(N\log N)$ complexity for problems from two dimensional
  scattering applications. We present numerical results for several
  test examples to illustrate the algorithm and its application to two
  dimensional high frequency scattering problems.
\end{abstract}



\section{Introduction}
\label{sec:intro}

\subsection{Problem statement}

In this paper, we consider the two dimensional $N$-body problem for
the high frequency Helmholtz kernel. Let $\{f_i, 1\le i \le N\}$ be a
set of charges located at points $\{p_i, 1\le i\le N\}$ in $\R^2$. We
assume that the points $\{p_i\}$ belong to a square centered at the
origin with size $K$. The problem is to evaluate the potentials
$\{u_i,1\le i \le N\}$ defined by
\begin{equation}
  u_i = \sum_{j=1}^N G(p_i,p_j) \cdot f_j
  \label{eq:nbody}
\end{equation}
where $G(x,y) = \frac{\i}{4} H^{(1)}_0(2\pi |x-y|)$ is the fundamental
solution of the 2D Helmholtz equation. In this paper, we use $\i$ to
denote $\sqrt{-1}$.

This computational problem mostly arises from the numerical solution
of 2D time harmonic scattering problems \cite{colton-1983-iemst}. For
example, suppose that $D \subset \R^2$ is a compact object with a
smooth boundary and $u^{inc}$ is the incoming field. If $D$ represents
a sound soft scatterer, the scattering field $u$ satisfies the
following Helmholtz equation with the Dirichlet boundary condition:
\[
-\Delta u - (2\pi)^2 u = 0 \quad\mbox{in}\; \R^d \setminus \bar{D}
\]
\[
u(x) = - u^{inc}(x) \quad\mbox{for}\; x\in\p D 
\]
\[
\lim_{r\rightarrow\infty} r \left( \frac{\p u}{\p r} - 2\pi\i u \right) = 0 
\]
where the wave number is set to be $2\pi$. The last condition is the
Sommerfeld radiation condition and guarantees that the scattering
field $u$ propagates to infinity. One highly efficient way to solve
this problem is to reformulate it into an equivalent boundary integral
equation (BIE)
\begin{equation}
  \frac{1}{2} \phi(x) + \int_{\p D} \left( \frac{\p G(x,y)}{\p n(y)} - \i \eta G(x,y) \right) \phi(y) d y
  = - u^{inc}(x)
  \label{eq:bie}
\end{equation}
where $n(y)$ is the exterior normal of $\p D$ at $y$, $\eta$ is some
fixed constant, and $\phi(x)$ for $x\in \p D$ is the unknown charge
distribution on the boundary $\p D$. Once $\phi$ is solved from
\eqref{eq:bie}, the scattering field $u$ can be simply computed with
an integral formula \cite{colton-1983-iemst}. The BIE approach has the
advantage of reducing the number of unknowns. The discrete version of
\eqref{eq:bie}, however, is a dense linear system which usually
requires an iterative algorithm like GMRES \cite{saad-1986-gmres} for
its solution. At each step of the iterative solver, we then need to
evaluate the computational problem in \eqref{eq:nbody}, with $\{p_i\}$
being the appropriate quadrature points.

It is well known that the complexity of a scattering problem often
scales with the size of scatterer in terms of the wavelength.  Since
the wavelength is taken to be 1 in our setup, the complexity of
\eqref{eq:nbody} depends on the number $K$, which can be of order
$10^4$ for a typical large scale scattering problem. Since one often
uses a constant number of points per wavelength when discretizing
\eqref{eq:bie}, the number of points $N$ is proportional to $K$.

\subsection{Previous work}

Direct computation of \eqref{eq:nbody} takes $O(N^2)$ operations. This
can be quite time consuming when $N$ is large. Various fast algorithms
have been proposed to reduce this complexity in the past two decades.
Among them, the most popular approach is the high frequency fast
multipole method (HF-FMM) developed by Rokhlin et al.
\cite{cheng-2006-riwfmm2d,rokhlin-1990-rsiest}. In the HF-FMM, the
whole computational domain is partitioned into a quadtree and one
associates with each square of the quadtree two expansions: the far
field expansion and the local field expansion
\cite{cheng-2006-riwfmm2d}.  These expansions allow one to accelerate
the computation in the low frequency region. In the high frequency
region, the Fourier transforms of these expansions are used instead to
achieve optimal efficiency since the translations between them become
diagonal operators under the Fourier basis. The HF-FMM has an optimal
$O(N\log N)$ complexity and has been widely used.

A different approach is to discrete the integral equation
\eqref{eq:bie} under the Galerkin framework using local Fourier bases
or wavelet packets. The stiffness matrix becomes approximately sparse
under these bases since most of the entries are close to zero and can
be safely discarded. Early algorithms
\cite{averbuch-2000-ecoi,bradie-1993-fnc,canning-1992-sasieok,deng-1999-fseie,deng-1999-cpwpb,golik-1998-wpfseie}
of this approach focus on finding the correct one dimensional basis,
while a recent development \cite{huybrechs-2006-twtmci} considers the
use of two dimensional wave packets which can offer more flexibility
and better compression rate.

Another early development is the multilevel matrix decomposition by
Michielssen and Boag \cite{michielssen-1996-mmda}. The three stage
multiplication algorithm, which is later named the butterfly algorithm
by \cite{oneil-2007-ncabft}, is quite similar to the FFT and brings
the overall complexity down to $O(N \log^2 N)$.

In \cite{engquist-2007-fdmaok}, we proposed an algorithm for the three
dimensional $N$-body problem of the high frequency Helmholtz kernel.
It relies on a low rank property of the 3D Helmholtz kernel for
certain geometric configurations. The algorithm organizes the
computation in a multidirectional and multilevel fashion and has an
optimal $O(N \log N)$ complexity.

\subsection{A multidirectional approach}

\begin{figure}
  \begin{center}
    \includegraphics[height=1.5in]{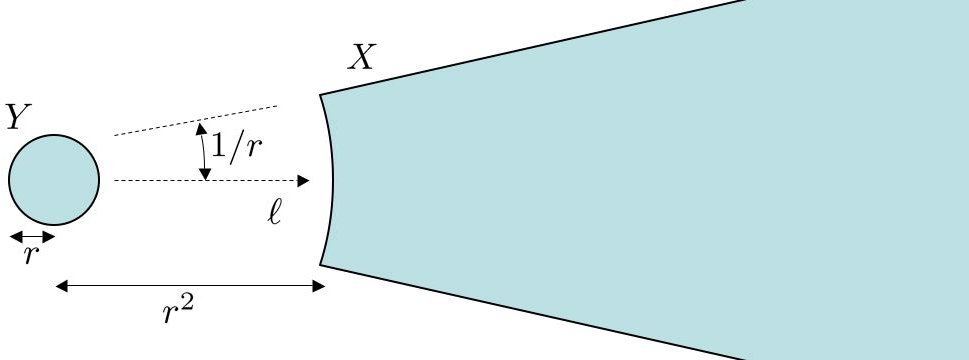}
  \end{center}
  \caption{Two sets $Y$ and $X$ that satisfy the directional parabolic
    separation condition.}
  \label{fig:xryr}
\end{figure}

In this paper, we adapt the approach in \cite{engquist-2007-fdmaok} to
the two dimensional $N$-body problem of the Helmholtz kernel. The main
idea is a similar low rank property of the 2D Helmholtz kernel. We say
that two sets $Y$ and $X$ satisfy the {\em directional parabolic
  separation condition} if $Y$ is a disk of radius $r$ and $X$ is the
set of points that belong to a cone with spanning angle $1/r$ and are
at least $r^2$ away from $Y$ (see Figure \ref{fig:xryr}).

Once $Y$ and $X$ satisfy the directional parabolic separation
condition, one can show that for any fixed accuracy the interaction
between $X$ and $Y$ via the Helmholtz kernel $G(x,y)$ is approximately
of low rank and the rank is independent of $r$. More precisely, for
any accuracy $\eps$, there exist a constant $T(\eps)$ and two sets of
functions $\{\alpha_i(x), 1\le i \le T(\eps) \}$ and $\{\beta_i(y),
1\le i \le T(\eps) \}$ such that for any $x\in X$ and $y\in Y$
\[
\left| G(x,y) - \sum_{i=1}^{T(\eps)} \alpha_i(x) \beta_i(y) \right| \le \eps
\]
(see Theorem \ref{thm:dlr}). Notice that $\{\alpha_i(x)\}$ and
$\{\beta_i(y)\}$ are only functions of $x$ and $y$ respectively. We
call such an approximation a {\em directional separated
  representation}. One major component of our approach is to use these
representations to build equivalent charges for well-separated
interaction. 

\begin{figure}
  \begin{center}
    \includegraphics[height=1.5in]{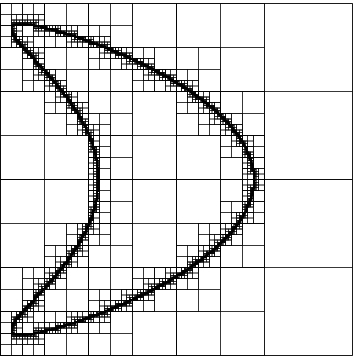} \hspace{0.1in}
    \includegraphics[height=1.5in]{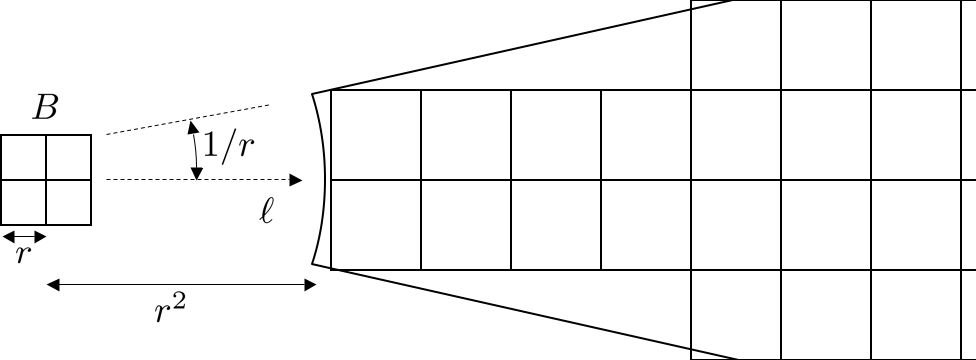}
  \end{center}
  \caption{Left: the quadtree constructed for a point distribution
    supported on a curve. Right: for each square $B$ in the high
    frequency region, its far field is partitioned into multiple
    wedges. We construct a low rank representation of the interaction
    between $B$ and each of its wedges. This representation is further
    used to accelerate the computation between $B$ and all the squares
    in the wedge.  }
  \label{fig:onewedge}
\end{figure}

Similar to the 3D algorithm in \cite{engquist-2007-fdmaok}, our 2D
algorithm starts by generating a quadtree for the whole computational
domain. In the low frequency region where the squares are of size less
than 1, the interactions are accelerated using the kernel independent
FMM algorithm in \cite{ying-2004-kiafmm}. In the high frequency region
where the squares are of size greater than or equal to 1, the far
field of each square is partitioned into wedges which follow the
directional parabolic separation condition (see Figure
\ref{fig:onewedge}). Between the square and each of its wedges, the
computation is accelerated via the directional separated
representation associated with the wedge.

Apart from extending the multidirectional algorithm of
\cite{engquist-2007-fdmaok} to the 2D Helmholtz kernel, this paper
also contains two new contributions:
\begin{itemize}
\item We provide an improved randomized procedure for the construction
  of the directional separated representations. The new procedure is
  more efficient and generates representations with smaller ranks.
\item Our algorithm has been applied to the solution of
  \eqref{eq:bie}.  This allows us to study large scatterers that are
  thousands of wavelengths wide.
\end{itemize}
  
The rest of this paper is organized as follows. In Section
\ref{sec:direct}, we briefly summarize the theoretical result on which
our approach is based and describe the new improved procedure for
constructing the separated representations. After describing our
algorithm for \eqref{eq:nbody} in detail in Section \ref{sec:algo}, we
present in Section \ref{sec:results} the numerical results for several
test examples. Finally, Section \ref{sec:concl} provides some comments
on future research directions. Though this paper focuses on the two
dimensional Helmholtz kernel, we would like to point out that our
algorithm is also applicable to other 2D oscillatory kernels such as
$e^{2\pi\i |x-y|}$.

\section{Directional Separated Representations}
\label{sec:direct}

\begin{definition}
  Let $f(x,y)$ be a function for $x \in X$ and $y \in Y$. We say
  $f(x,y)$ has a $T$-term $\eps$-expansion for $X$ and $Y$ if there
  exist functions $\{ \alpha_i(x), 1 \le i \le T\}$ and $\{
  \beta_i(y), 1 \le i \le T\}$ such that
  \[
  \left| f(x,y) - \sum_{i=1}^{T} \alpha_i(x) \beta_i(y) \right| \le \eps
  \]
  for all $x \in X$ and $y \in Y$.
\end{definition}
Since the two sets of functions $\{\alpha_i(x)\}$ and $\{\beta_i(y)\}$
depend only on $x$ and $y$ respectively, the above expansion is called
{\em separated}. Suppose $r \ge \sqrt{2}$. For our problem, we take
\begin{equation}
  Y = B(0,r)
  \quad\mbox{and}\quad
  X = \{ x: \theta(x,\ell) \le 1/r, |x| \ge r^2\}
  \label{eq:xryr}
\end{equation}
where $\ell$ is a given unit vector and $\theta(a,b)$ is the spanning
angle between vectors $a$ and $b$. The geometric relationship between
$Y$ and $X$ is illustrated in Figure \ref{fig:xryr}. The following
theorem serves as the theoretical foundation of our approach.

\begin{theorem}
  For any $\eps>0$, there exists a number $T(\eps)$ which is
  independent of $r$ such that 
  \[
  G(x,y) = \frac{\i}{4} H^{(1)}_0(2\pi |x-y|)
  \]
  has a $T(\eps)$-term $\eps$-expansion for any $X$ and $Y$ given
  by \eqref{eq:xryr}.
  \label{thm:dlr}
\end{theorem}

The representation guaranteed by Theorem \ref{thm:dlr} is called a
{\em directional separated representations} for the obvious reason.
One way to prove this theorem is to use the asymptotic behavior of
$H^{(1)}_0$ for large arguments
\cite{abramowitz-1992-hmf,bronshtein-1997-hm}:
\[
H^{(1)}_0 (r) = \sqrt{ \frac{2}{\pi r} }  \left( e^{\i (r-\pi/4)} + O\left(\frac{1}{r}\right) \right),
\]
and then follow the same path as the proof for Theorem 2.2 in
\cite{engquist-2007-fdmaok}.

\subsection{Construction of directional separated representation}

\begin{figure}
  \begin{center}
    \includegraphics[height=1.5in]{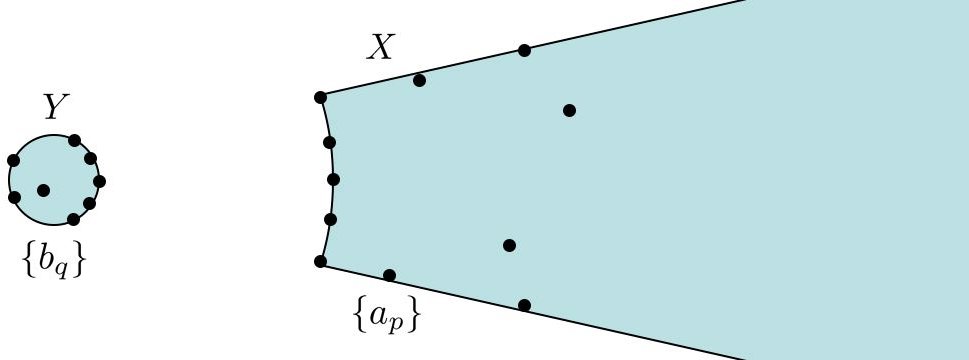}
  \end{center}
  \caption{Constructions of the separated representation between $X$ and $Y$. 
    $\{b_q\}$ are the samples associated with the columns in $A_c$ (Step 2).
    $\{a_p\}$ are the samples associated with the columns in $A_r$ (Step 3).
  }
  \label{fig:samples}
\end{figure}

A procedure based on random sampling has been described in
\cite{engquist-2007-fdmaok} for the construction of these directional
separated representations. In the rest of this section, we propose an
improved version which gives lower separation ranks and better
accuracy based on our numerical experience. For a given pair $Y$ and
$X$ that satisfy the directional parabolic separation condition, our
new procedure takes the following steps:

\begin{enumerate}
\item Sample $Y$ randomly with a set of samples $\{y_j, 1\le j \le
  N_Y\}$. In our implementation, we use 2 to 3 points per wavelength
  and the number of samples $N_Y$ grows linearly with the area of $Y$.
  Sample $X$ similarly with a set of samples $\{x_i, 1\le i \le
  N_X\}$. Let $A$ be the matrix defined by
  \[
  A_{ij} = G(x_i, y_j) = \frac{\i}{4} H^{(1)}_0(2\pi |x_i-y_j|),
  \]
  for $ 1\le i \le N_X$ and $1\le j \le N_Y$. In the language of
  linear algebra, Theorem \ref{thm:dlr} states that $A$ can be
  factorized, within error $O(\eps)$, into the product of two
  matrices, the first containing $T(\eps)$ columns and the second
  containing $T(\eps)$ rows. 
  
\item Let $A_1$ be the submatrix of $A$ containing a set of $N_1$
  randomly selected rows. Here we set $N_1 \approx 3 \cdot T(\eps)$ in
  practice.  Our goal is to find a set of $T(\eps)$ columns of $A_1$
  that has the largest $T(\eps)$-dimensional volume. Since $A_1$ is
  only of size $O(T(\eps)) \times N_Y$, one can use either the
  interpolative decomposition \cite{cheng-2005-clrm} or the greedy
  standard pivoted QR factorization to find these columns. Both
  algorithms have an $O(N_Y)$ complexity. Suppose the pivoted QR
  factorization is used.  We then have the decomposition
  \[
  A_1 P_1 = Q_1 R_1,
  \]
  where $P_1$ is a permutation matrix, $Q_1$ is orthonormal, and $R_1$
  is upper triangular. Now identify the diagonal elements of $R_1$
  which are less than $\eps$ and truncate the associated columns of
  $Q_1$ and rows of $R_1$. Denote the resulting matrices by $Q_{1,c}$
  and $R_{1,c}$. Since $A_1$ itself has an $O(T(\eps))$-expansion,
  $Q_{1,c}$ contains only $O(T(\eps))$ columns in practice. Moreover,
  it is clear that
  \[
  Q_{1,c} R_{1,c} = A_{1,c},
  \]
  where $A_{1,c}$ is the submatrix containing the columns of $A_1$
  from which the matrix $Q_{1,c}$ is generated. We denote by $A_c$ the
  submatrix of $A$ that consists of the same columns. The $O(T(\eps))$
  samples of $Y$ associated with these columns are denoted $\{b_q\}$.
  
\item Let $A_2$ be a submatrix of $A$ containing a set of $N_2$
  randomly selected columns. We again set $N_2 \approx 3 \cdot T(\eps)$.
  Repeat the previous step on $A_2^*$.  As a result, we have two
  matrices $Q_{2,r}$ and $R_{2,r}$. $Q_{2,r}$ is orthonormal and has
  $O(T(\eps))$ columns again, while $R_{2,r}$ is upper triangular.
  They satisfy the relationship
  \[
  R_{2,r}^* Q_{2,r}^* = A_{2,r},
  \]
  where $A_{2,r}$ is a submatrix containing appropriate rows of $A$.
  We denote by $A_r$ the submatrix of $A$ that consists of the same
  rows and by $\{a_p\}$ the $O(T(\eps))$ samples of $X$ associated
  with these rows (see Figure \ref{fig:samples}).
  
\item We randomly pick a set $S$ of $N_S$ rows and a set $T$ of $N_T$
  columns. In practice, we choose $N_S$ and $N_T$ to be equal to
  $10 \cdot T(\eps)$. Set $A_3$ to be the minor containing the elements from
  rows in $S$ and columns in $T$, $A_{c,S}$ to be the submatrix of
  $A_c$ containing the rows in $S$, and $A_{r,T}$ to be the submatrix
  of ${A_r}$ containing the columns in $T$. Next, we choose $ D =
  (A_{c,S})^+ A_3 (A_{r,T})^+, $ where $(\;)^+$ stands for
  pseudoinverse. We claim that
  \[
  \left| A - A_c D A_r \right| = O(\eps).
  \]
  Such an approximate factorization is often called a pseudoskeleton 
  approximation of $A$ in the literature (see
  \cite{goreinov-1997-tpa,goreinov-1997-pasgs}). Notice that the
  matrix $D$ has only $O(T(\eps))$ rows and columns. Denoting the
  entries of $D$ by $d_{qp}$, we can rewrite the previous statement in
  the form
  \[
  \left|
    G(x_i,y_j)    - \sum_{p,q}     G(x_i,b_q)    \cdot  d_{qp}\cdot    G(a_p,y_j)
  \right| = O(\eps)
  \]
  for all $x_i$ and $y_j$.

\item Finally, since $\{x_i\}$ and $\{y_j\}$ sample the sets $X$ and
  $Y$ with a constant number of points per wavelength, it is
  reasonable to expect
  \begin{equation}
  \left|
    G(x,y)    - \sum_{p,q}     G(x,b_q) \cdot     d_{qp}\cdot    G(a_p,y)
  \right| = O(\eps)
  \label{eq:fnlexp}
  \end{equation}
  for any $x \in X \cap B(0,K)$ and $y \in Y$.
\end{enumerate}

Since both $\{a_p\}$ and $\{b_q\}$ are of order $O(T(\eps))$, it is
clear that \eqref{eq:fnlexp} is a low rank separated representation.
Moreover, we only need to store $\{a_p\}$, $\{b_q\}$, and $D$ for
\eqref{eq:fnlexp}, thus reducing the storage requirement dramatically.
We would like to point out that recently there has been a lot of
research devoted to problems similar to \eqref{eq:fnlexp} (see
\cite{bebendorf-2003-alacm,drineas-2006-fmcam2,drineas-2006-fmcam3,martinsson-2006-rafam}
for details).

This randomized procedure performs quite well in practice as we will
see from the numerical results in Section \ref{sec:results}. Though we
do not yet have a proof, the following heuristic argument provides
some useful insights. In the standard pseudoskeleton approximation
\cite{goreinov-1997-tpa,goreinov-1997-pasgs}, an $m \times n$ matrix
$A$ has the following approximation:
\[
A \approx A_c G A_r,
\]
where $A_c$, $G$, and $A_r$ are of size $m \times k$, $k \times k$,
and $k \times n$ respectively. Often $A_c$ contains the columns of $A$
that have the largest $k$-dimensional volume and, similarly, $A_r$
contains the rows with the largest $k$-dimensional volume. Finding
these columns and rows are quite expensive if both $m$ and $n$ are
large. Suppose now that we can project the columns (or rows) of $A$
onto a $p$ dimensional subspace $L$ which is randomly selected from
all $p$-dimensional subspaces with the uniform rotational invariant
probability measure. As long as $p$ is adequately larger than $k$, the
volume spanned by any set of $k$ columns (or rows) is preserved to a
good accuracy \cite{dasgupta-2003-eptjl,magen-2002-dr}. Therefore, one
efficient method to find the columns of $A$ with the largest volume
would be to
\begin{enumerate}
\item project $A$ onto a random $p$ dimensional subspace, 
\item find the columns of the projected matrix that have the largest
  $k$-dimensional volume,
\item pick the corresponding columns of $A$ to be the answer.
\end{enumerate}
The only difference between this approach and the second and third
steps of our randomized procedure is that we only project to a random
set of coordinates, which is much more restrictive than the uniform
random projection. However, since both the columns and the rows of our
matrix $A$ is highly oscillatory and incoherent with the Dirac
functions, our procedure works well in practice.

\subsection{Equivalent charges}

The directional separated representation \eqref{eq:fnlexp} provides a
way to represent the potential in $X$ generated by the charges inside
$Y$ in a compact way. Suppose that $X$ is centered around the unit
direction $\ell$ and $\{f_i\}$ are the charges located at points
$\{y_i\}$ in $Y$. After applying \eqref{eq:fnlexp} to $y=y_i$ for each
$y_j$ and summing them up with weight $f_i$, we have
\[
\left|
  \sum_i G(x,y_i) f_i - 
  \sum_q G(x,b_q)
  \left(
  \sum_p d_{qp}
  \sum_i G(a_p,y_i) f_i 
  \right)
\right| = O(\eps).
\]
This states that we can place a set of charges 
\begin{equation}
  \left\{
    \sum_p d_{qp}
    \sum_i G(a_p,y_i) f_i
  \right\}
  \label{eq:doed}
\end{equation}
at points $\{b_q\}$ in order to reproduce the potential generated by
the charges $\{f_i\}$ located at points $\{y_i\}$. We call the charges
in \eqref{eq:doed} the {\em directional outgoing equivalent charges}
of $Y$ in direction $\ell$ and the points $\{b_q\}$ the {\em
  directional outgoing equivalent points} of $Y$ in direction $\ell$.
In addition, we refer to the quantities
\begin{equation}
  \left\{
    \sum_i G(a_p,y_i) f_i
  \right\}
  \label{eq:docp}
\end{equation}
as the {\em directional outgoing check potentials} of $Y$ in direction
$\ell$ and the points $\{a_p\}$ as the {\em directional outgoing check
  points} of $Y$ in direction $\ell$. Given the check potentials, the
equivalent charges can be computed easily by a multiplication with
$D$.

Let us now reverse the role of $X$ and $Y$. Suppose we have a set of
charges $\{f_i\}$ located at points $\{x_i\}$ in $X$. Since $G(x,y) =
G(y,x)$,
\[
\left|
  \sum_i G(y,x_i) f_i - 
  \sum_p G(y,a_p)
  \sum_q d_{qp}
  \sum_i G(b_q,x_i) f_i
\right| = O(\eps).
\]
This states that we can put a set of charges 
\begin{equation}
  \left\{
    \sum_q d_{qp}
    \sum_i G(b_q,x_i) f_i 
  \right\}
  \label{eq:died}
\end{equation}
at points $\{a_p\}$ and they reproduce the potential generated by the
charges $\{f_i\}$ located at points $\{x_i\}$. Therefore, we call the
charges in \eqref{eq:died} the {\em directional incoming equivalent
  charges} of $Y$ in direction $\ell$ and the locations $\{a_p\}$ the
{\em directional incoming equivalent points} of $Y$ in direction
$\ell$. In analogy to the previous terminology,
\begin{equation}
  \left\{
    \sum_i G(b_q,x_i) f_i
  \right\}
  \label{eq:dicp}
\end{equation}
are called the {\em directional incoming check potentials} of $Y$ in
direction $\ell$ and the location $\{b_q\}$ are called the {\em
  directional incoming check points} of $Y$ in direction $\ell$.

\section{Algorithm Description}
\label{sec:algo}


Without loss of generality, we assume
that the size of the domain $K = 2^{2L}$ for a positive integer $L$.

\subsection{Data structure}
We start by constructing a quadtree which contains the whole
computational domain. We often use $B$ to denote a square in the
quadtree and $w$ for its width. A square $B$ is said to be in the low
frequency regime if $w < 1$ and in the high frequency regime if $w \ge
1$. In the high frequency regime of the quadtree, no adaptivity is
used, i.e., every non-empty square is further partitioned until the
width of the square is less than $1$. In the low frequency regime, a
square $B$ is partitioned as long as the number of points in $B$ is
greater than a fixed constant $N_p$. The value of $N_p$ is chosen to
optimize the computational complexity and, in practice, we pick $N_p =
50$.

For a square $B$ in the low frequency regime, its data structure follows
the description of the kernel independent FMM in
\cite{ying-2004-kiafmm}. The near field $N^B$ is the union of the
squares $A$ that satisfies $dist(A,B) = 0$, where $dist(A,B) =
\inf_{x\in A,y\in B} |x-y|$. The far field $F^B$ is the complement of
$N^B$. The interaction list $I^B$ contains all the squares in $N^P
\backslash N^B$ on $B$'s level, where $P$ is the parent square of $B$.
\begin{itemize}
\item $\{y^{B,o}_k\}$, $\{f^{B,o}_k\}$, $\{x^{B,o}_k\}$ and
  $\{u^{B,o}_k\}$ are, respectively, the {\em outgoing} equivalent
  points, equivalent charges, check points, and check potentials.
\item $\{y^{B,i}_k\}$, $\{f^{B,i}_k\}$, $\{x^{B,i}_k\}$ and
  $\{u^{B,i}_k\}$ are, respectively, the {\em incoming} equivalent
  points, equivalent charges, check points, and check potentials.
\end{itemize}

\begin{figure}
  \begin{center}
    \includegraphics[height=1.8in]{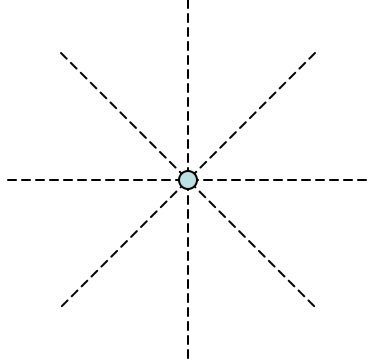}
    \includegraphics[height=1.8in]{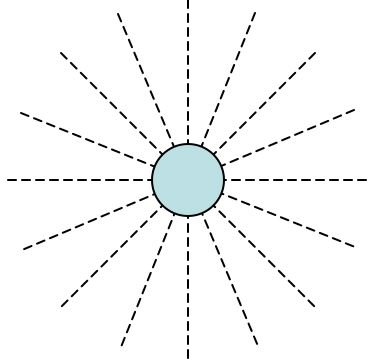}
    \includegraphics[height=1.8in]{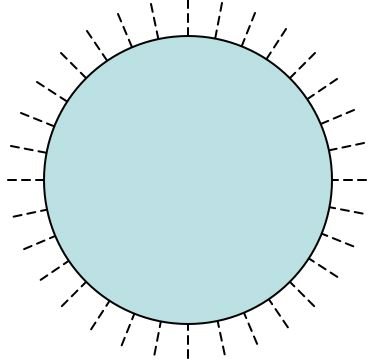}
  \end{center}
  \caption{The far field is partitioned into wedges. From left to
    right, $w=1,2,4$. The radii are 1,4, and 16, respectively.
  }
  \label{fig:six}
\end{figure}

For a square $B$ in the high frequency region, the near field $N^B$ is
the union of all the squares $\{A\}$ that satisfy $dist(A,B) \le w^2$.
The far field $F^B$ is the complement of $N^B$. The interaction list
$I^B$ contains all the squares in $N^P \backslash N^B$ on $B$'s level,
where $P$ is $B$'s parent square. Notice that the far field of a
square $B$ in the high frequency region is pushed away in order to be
compatible with the directional parabolic separation condition. The
far field $F^B$ is further partitioned into a group of directional
wedges, each belonging to a cone with spanning angle $O(1/w)$. We
denote the set of all the wedges of $B$ by $\{W^{B,\ell}\}$.  In
Figure \ref{fig:six}, we illustrate the case for for $w=1,2,4$.

For each square $B$ and each direction $\ell$, we summarize the relevant
quantities as follows:
\begin{itemize}
\item $\{y^{B,o,\ell}_k\}$, $\{f^{B,o,\ell}_k\}$,
  $\{x^{B,o,\ell}_k\}$, and $\{u^{B,o,\ell}_k\}$ are the {\em outgoing
    directional} equivalent points, equivalent charges, check points
  and check potentials respectively.
\item $\{y^{B,i,\ell}_k\}$, $\{f^{B,i,\ell}_k\}$,
  $\{x^{B,i,\ell}_k\}$, and $\{u^{B,i,\ell}_k\}$ are the {\em incoming
    directional} equivalent points, equivalent charges, check points
  and check potentials respectively.
\end{itemize}

\subsection{Translation operators}

Following the convention in
\cite{greengard-1987-afaps,rokhlin-1990-rsiest}, we name these
operators M2M, L2L, and L2L translations, though no multipole or local
expansions are involved in our algorithm. The translation operators
for squares in the low frequency regime are detailed already in
\cite{ying-2004-kiafmm}. The operators in the high frequency regime
are more complicated. The main reason is that the computations are now
directional.

For a square $B$ in the high frequency regime, the {\em M2M translation}
constructs the outgoing directional equivalent charges of $B$ from the
outgoing equivalent charges of $B$'s children. There are two cases to
consider. In the first case, $w=1$.  The children squares have only
nondirectional equivalent charges. The M2M translation iterates over
all of the directional indices $\{\ell\}$ of $B$, and the steps for a
fixed direction $\ell$ are as follows:
\begin{enumerate}
\item Use $\bigcup_{C} \{y^{C,o}_k\}$ as source points in $B$ and
  $\bigcup_{C} \{f^{C,o}_k\}$ as source charges. Here the union is
  taken over all of the children squares of $B$.
\item Compute $\{u^{B,o,\ell}_k\}$ at points $\{x^{B,o,\ell}_k\}$ with
  kernel evaluation, and then obtain $\{f^{B,o,\ell}_k\}$ by
  multiplying $\{u^{B,o,\ell}_k\}$ with the matrix $D$ associated with
  the wedge $W^{B,\ell}$.
\end{enumerate}

\begin{figure}
  \begin{center}
    \includegraphics[height=1.5in]{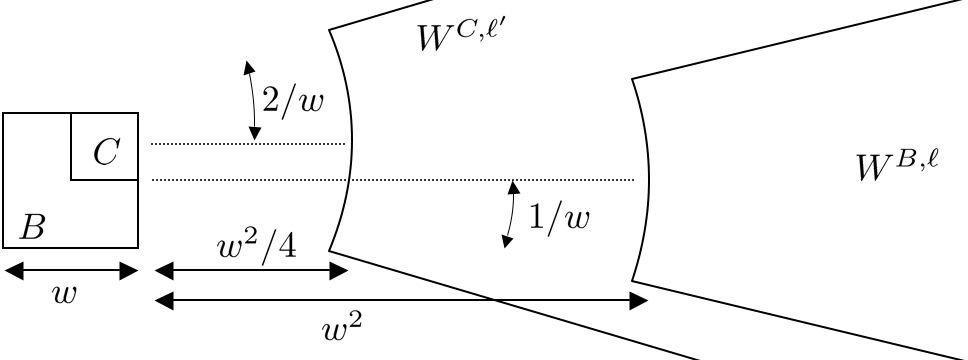}
  \end{center}
  \caption{$B$ is a square with width $w>1$. For any fixed $\ell$,
    there exists $\ell'$ such that $W^{B,\ell}$ is contained in
    $W^{C,\ell'}$ where $C$ is any one of $B$'s children.}
  \label{fig:ellell}
\end{figure}

In the second case, $w>1$. Now the children squares have directional
equivalent charges as well. The M2M translation iterates over all of
the directional indices $\{\ell\}$ of $B$. The steps for a fixed
direction $\ell$ are as follows:
\begin{enumerate}
\item Pick $\ell'$, a direction associated with the squares of width
  $w/2$, such that the wedge $W^{B,\ell}$ is contained in the wedge
  $W^{C,\ell'}$ where $C$ stands for anyone of $B$'s children.  The
  existence of $\ell'$ is ensured by the way we partition $F^B$
  (see Figure \ref{fig:ellell}).
\item Use $\bigcup_{C} \{ y^{C,o,\ell'}_k\}$ as source points in $B$
  and $\bigcup_{C} \{f^{C,o,\ell'}_k\}$ as source charges. Here the
  union is taken over all the children squares of $B$.
\item Compute $\{u^{B,o,\ell}_k\}$ at $\{x^{B,o,\ell}_k\}$ with kernel
  evaluation and then obtain $\{f^{B,o,\ell}_k\}$ by multiplying
  $\{u^{B,o,\ell}_k\}$ with the matrix $D$ associated with the wedge
  $W^{B,\ell}$.
\end{enumerate}

The {\em L2L translation} constructs the incoming check potentials of
$B$'s children from the incoming directional check potentials of $B$.
Again there are two cases to consider. In the first case $w=1$. The
children squares have only nondirectional check potentials. The L2L
translation iterates over all of the directional indices $\{\ell\}$ of
$B$, and the steps for a fixed direction $\ell$ are as follows:
\begin{enumerate}
\item Compute $\{f^{B,i,\ell}_k\}$ from $\{u^{B,i,\ell}_k\}$ by multiplying
  it with the appropriate $D$ matrix.
\item For each child $C$ of the square $B$, add to $\{u^{C,i}_k\}$ the
  potentials evaluated at $\{x^{C,i}_k\}$ using $\{f^{B,i,\ell}_k\}$
  as the source charges at $\{y^{B,i,\ell}_k\}$.
\end{enumerate}

In the second case, $w>1$. Now the children squares have directional
equivalent charges. The L2L translation iterates over all of the
directional indices $\{\ell\}$ of $B$. The steps for a fixed direction
$\ell$ are as follows:
\begin{enumerate}
\item Pick $\ell'$, a direction associated with the squares of width
  $w/2$, such that the wedge $W^{B,\ell}$ is contained in the wedge
  $W^{C,\ell'}$ where $C$ stands for anyone of $B$'s children. 
\item Compute $\{f^{B,i,\ell}_k\}$ from $\{u^{B,i,\ell}_k\}$ by
  multiplying it with the appropriate $D$ matrix.
\item For each child $C$ of the square $B$, add to $\{u^{C,i,\ell'}_k\}$
  the potentials evaluated at $\{x^{C,i,\ell'}_k\}$ using
  $\{f^{B,i,\ell}_k\}$ as the source charges at
  $\{y^{B,i,\ell}_k\}$.
\end{enumerate}

Finally, the {\em M2L translation} is applied to pairs of squares $A$
and $B$ on the same level of the quadtree. They need to be on each
other's interaction lists. Suppose $B$ falls into the wedge
$W^{A,\ell}$ of $A$ while $A$ falls into the wedge $W^{B,\ell'}$ of
$B$. The implementation of the M2L translation contains only one step:
\begin{enumerate}
\item Add to $\{u^{B,i,\ell'}_k\}$ the potentials evaluated at
  $\{x^{B,i,\ell'}_k\}$ using the charges $\{f^{A,o,\ell}_k\}$
  at points $\{y^{A,o,\ell}_k\}$.
\end{enumerate}

To summarize the discussion on the transition operators, we would like
to emphasize that all of these operators involve only kernel evaluation
and matrix-vector multiplication with precomputed matrices.
Therefore, they are simple to implement and highly efficient.

\subsection{Algorithm}
\label{sec:algo-algo}

\begin{figure}
  \begin{center}
    \includegraphics[height=3.2in]{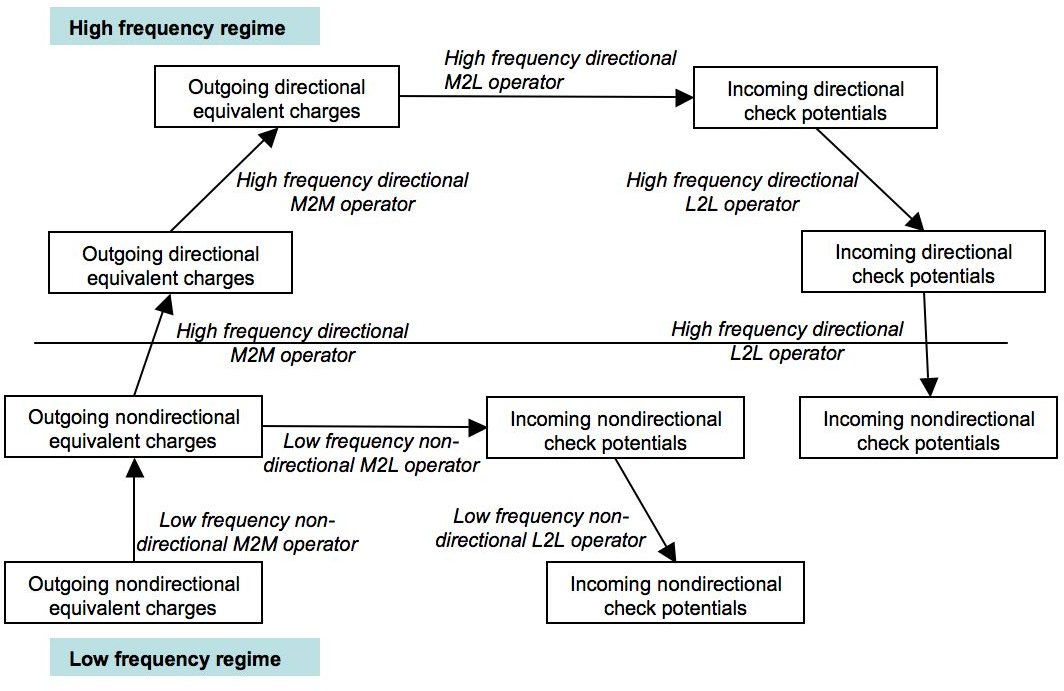}
  \end{center}
  \caption{A small part of the quadtree used in the computation. Each
    rectangular region stands for a square of the quadtree. The diagram
    shows how the outgoing nondirectional equivalent charges from a
    leaf square have been transformed into incoming nondirectional check
    potentials at other leaf squares. Far field interaction involves
    directional computation in the high frequency regime.}
  \label{fig:algo}
\end{figure}

Now we are ready to give the overall structure of our new algorithm.
It has exactly the same structure as the 3D algorithm in
\cite{engquist-2007-fdmaok} and we simply reproduce it here:
\begin{enumerate}
\item Construct the quadtree. In the high frequency regime, the squares are
  partitioned uniformly. In the low frequency regime, a leaf square
  contains at most $N_p$ points.
\item Travel up in the quadtree and visit the squares in the low frequency
  regime. These squares have width less than 1. For each square $B$,
  compute its outgoing nondirectional equivalent charges
  $\{f^{B,o}_k\}$. This is done using the low frequency
  nondirectional M2M translation.
\item Travel up in the quadtree and visit the squares in the high
  frequency regime. For every such square $B$, use the high frequency
  directional M2M translation to compute the outgoing directional
  equivalent charges $\{f^{B,o,\ell}_k\}$ for each outgoing direction
  $\ell$. We skip the squares with width greater than $\sqrt{K}$ since
  their interaction lists are empty.
\item Travel down in the quadtree and visit the squares in the high
  frequency regime. For every such square $B$ and for each direction
  $\ell$, perform the following two steps:
  \begin{enumerate}
  \item Transform the outgoing directional equivalent charges
    $\{f^{A,o,\ell}_k\}$ of all of the squares $\{A\}$ in $B$'s
    interaction list and in direction $\ell$ via the high frequency
    directional M2L translation. Next, add the result to the incoming
    directional check potentials $\{u^{B,i,\ell}_k\}$.
  \item Perform the high-frequency directional L2L translation to
    transform $\{u^{B,i,\ell}_k\}$ to the incoming check potentials
    for $B$'s children.
  \end{enumerate}
  Again, we skip the squares with width greater than $\sqrt{K}$.
\item Travel down in the quadtree. For every square $B$ in the low
  frequency regime, we perform the following two steps:
  \begin{enumerate}
  \item Transform the outgoing nondirectional equivalent charges
    $\{f^{A,o}_k\}$ of all of the squares $\{A\}$ in $B$'s interaction
    list via the low frequency nondirectional M2L operator. Next, add
    the result to the incoming nondirectional check potentials
    $\{u^{B,i}_k\}$.
  \item Perform the low frequency directional L2L translation.
    Depending on whether $B$ is a leaf square or not, add the result to
    the incoming check potentials of $B$'s children or to the
    potentials at the original points inside $B$.
  \end{enumerate}
\end{enumerate}

An illustration of the various components of the algorithm is given in
Figure \ref{fig:algo}. The following theorem summarizes the
complexity of the proposed algorithm.

\begin{theorem}
  Let $\S$ be a rectifiable curve in $B(0,1/2)$. Suppose that for a
  fixed $K$ the points $\{p_i, 1 \le i \le N\}$ are samples of $K\S$,
  where $N=O(K)$ and $K\S = \{K \cdot p, p\in \S\}$ (the surface
  obtained by magnifying $\S$ by a factor of $K$). Then, for any
  prescribed accuracy, the proposed algorithm has a computational
  complexity $O(K \log K) = O(N\log N)$.
\end{theorem}

The proof of this theorem follows closely the steps of Theorem 4.1 of
\cite{engquist-2007-fdmaok}. The main step of the proof is the
observation that, for any fixed $w>1$, there are at most $O(K/w)$
squares of size $w$ and, for each of them, there are at most $O(w)$
squares for which we apply the M2L operator.

\section{Numerical Results}
\label{sec:results}

In this section, we provide some numerical results to illustrate the
properties of our new algorithm. All of the computational results
below are obtained on a desktop computer with a 2.8 GHz CPU.

Let us first study the performance of the randomized procedure
presented in Section \ref{sec:direct}. In Table \ref{tbl:rank}, we
list the number of terms in the separated representation for two sets
$X$ and $Y$ for different choices of accuracy $\eps$ and square width
$w$. Here $r$, the radius of $Y$, is set to be $\sqrt{2} w$ so that
the square of width $w$ is contained in $Y$. We can see from Table
\ref{tbl:rank} that the separation rank is bounded by a constant which
is independent of the values of $w$. This is consistent with our
theoretical estimate in Theorem \ref{thm:dlr}.  In fact, as $w$ grows,
it seems that the separation rank decays slightly.

\begin{table}
  \begin{center}
    \begin{tabular}{|c|cccccccc|}
      \hline
      & $w=1$ & $w=2$ & $w=4$ & $w=8$ & $w=16$ & $w=32$ & $w=64$ & $w=128$ \\
      \hline
      $\eps$=1e-4 & 14 & 11 & 11 & 10 & 9 & 9 & 9 & 9\\
      $\eps$=1e-6 & 19 & 16 & 14 & 13 & 12 & 12 & 12 & 11\\
      $\eps$=1e-8 & 27 & 20 & 16 & 15 & 15 & 15 & 14 & 14\\
      \hline
    \end{tabular}
  \end{center}
  \caption{The separation rank of the directional separated representation
    for different choices of requested accuracy $\eps$ and square size $w$.}
  \label{tbl:rank}
\end{table}

Next, we applied our algorithm to the $N$-body problems on several
objects. In our experiments, the boundary of each object is
represented by a piecewise smooth curve.  For these tests, the point
set $\{p_i\}$ is generated by sampling the curve randomly with about
$20$ points per wavelength. The densities $\{f_i\}$ are generated from
a random distribution with mean $0$. We use $\{u_i\}$ to denote the
true discrete potentials and $\{u_i^a\}$ to denote the approximations
obtained through our algorithm. We estimate the relative error by
picking a set $S$ of $200$ points from $\{p_i\}$.  The true potentials
$\{u_i, i\in S\}$ are computed by using direct evaluation.  The error
is then estimated to be
\[
\sqrt{
  \frac{ \sum_{i\in S} |u_i - u_i^a |^2 } { \sum_{i\in S} |u_i|^2 }
}.
\]

Before reporting the results, let us summarize the notations we use
here: $N$ is the number of points, $K$ is the size of the problem in
terms of the wavelength, $\eps$ is the prescribed error threshold such
that the final error is to be bounded by a constant multiple of
$\eps$, $T_a$ is the running time of our algorithm in seconds, $T_d$
is the running time of the direct evaluation in seconds, $T_d / T_a$
is the speedup factor, and $\eps_a$ is the resulting error of our
algorithm.

The first example is a circle and the results are
summarized in Table \ref{tbl:circ}. The second example is an airfoil
and the results are shown in Table \ref{tbl:foil}. The final example
is a kite-shaped object and we report the numbers in Table
\ref{tbl:kite}. These numbers demonstrate clearly that our algorithm
scales exactly like $O(N\log N)$ in terms of the number of points.
Furthermore, the error seems to grow only slightly as we increase the
number of points.

\begin{table}[h]
  \begin{center}
    \includegraphics[height=2in]{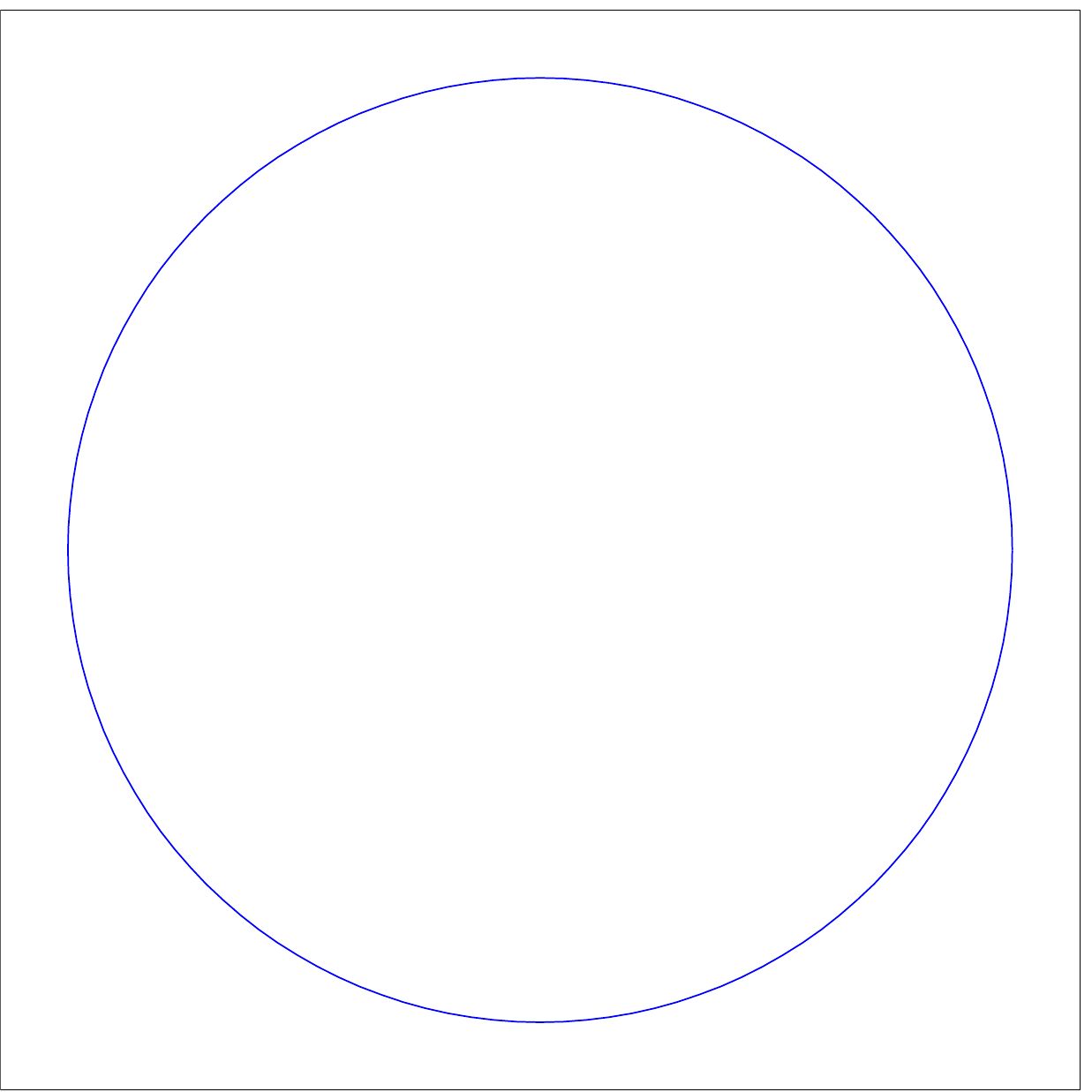}
    \includegraphics[height=2in]{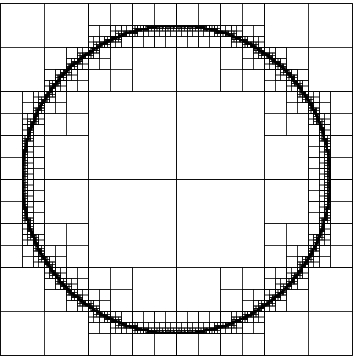}\\
    \vspace{0.1in}
    \begin{tabular}{|c|ccccc|}
      \hline
      $(K,\eps)$ & $N$ & $T_a$(sec) & $T_d$(sec) & $T_d / T_a$ & $\eps_a$\\
      \hline
      (2048,1e-4) & 1.13e+5 & 3.40e+1 & 8.05e+3 & 2.37e+2 & 1.25e-4\\
      (8192,1e-4) & 4.50e+5 & 1.56e+2 & 1.28e+5 & 8.21e+2 & 1.31e-4\\
      (32768,1e-4)& 1.80e+6 & 7.07e+2 & 2.06e+6 & 2.91e+3 & 1.80e-4\\
      \hline
      (2048,1e-6) & 1.13e+5 & 5.30e+1 & 8.00e+3 & 1.51e+2 & 7.88e-7\\
      (8192,1e-6) & 4.50e+5 & 2.39e+2 & 1.28e+5 & 5.37e+2 & 9.98e-7\\
      (32768,1e-6)& 1.80e+6 & 1.08e+3 & 2.06e+6 & 1.91e+3 & 1.00e-6\\
      \hline
      (2048,1e-8) & 1.13e+5 & 8.20e+1 & 8.05e+3 & 9.82e+1 & 8.48e-9\\
      (8192,1e-8) & 4.50e+5 & 3.57e+2 & 1.29e+5 & 3.60e+2 & 1.18e-8\\
      (32768,1e-8)& 1.80e+6 & 1.58e+3 & 2.07e+6 & 1.31e+3 & 1.30e-8\\
      \hline
    \end{tabular}
  \end{center}
  \caption{Results of a circle with the Helmholtz kernel.
    $N$ is the number of points, $K$ is the size of the problem
    in terms of the wavelength, $\eps$ is the prescribed error threshold
    such that the final error is to be bounded by a constant multiple of
    $\eps$, $T_a$ is the running time of our algorithm in seconds, $T_d$
    is the running time of the direct evaluation in seconds, $T_d / T_a$
    is the speedup factor, and $\eps_a$ is the estimated error of our
    algorithm.}
  \label{tbl:circ}
\end{table}

\begin{table}[h]
  \begin{center}
    \includegraphics[height=2in]{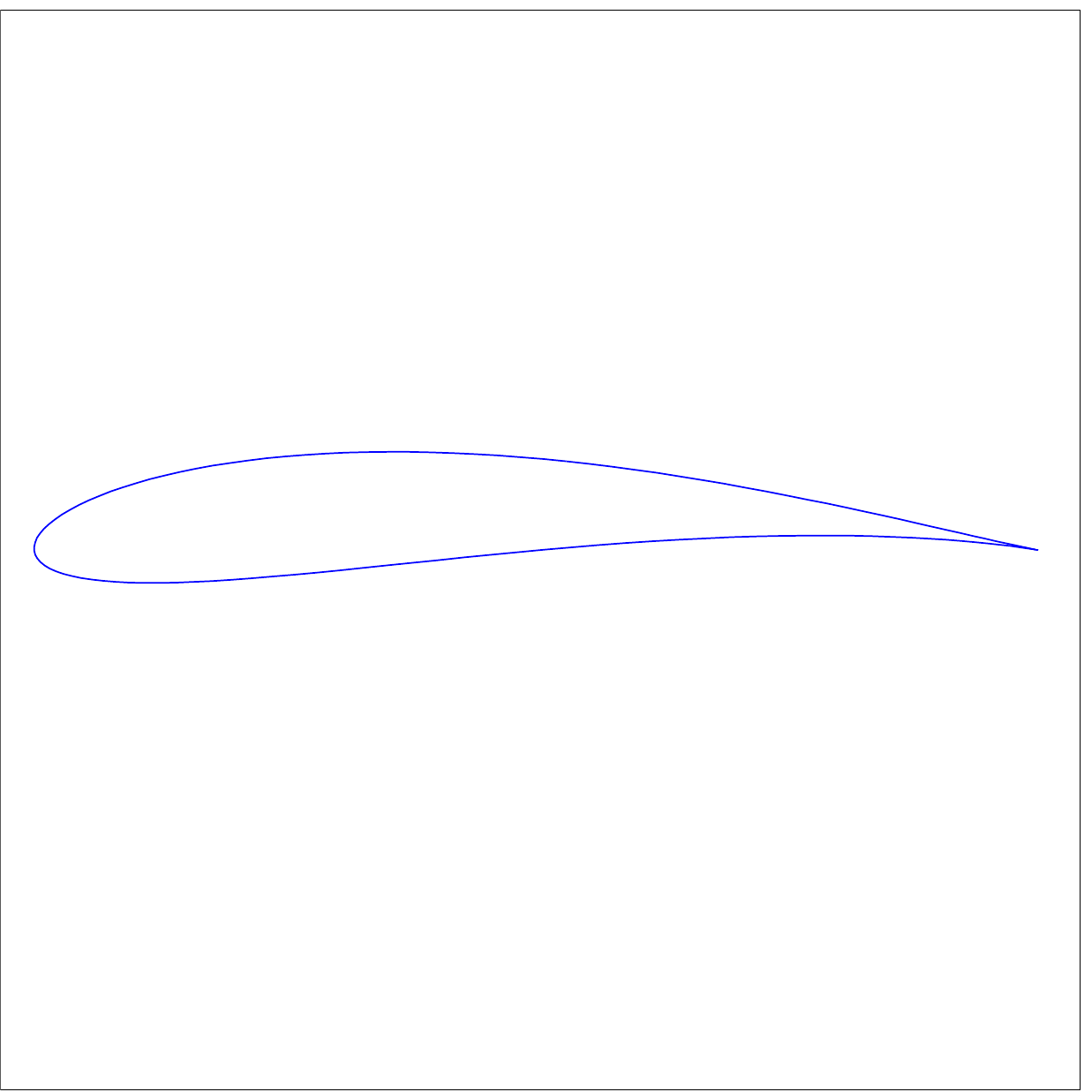}
    \includegraphics[height=2in]{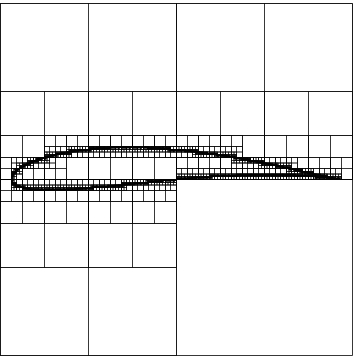}\\
    \vspace{0.1in}
    \begin{tabular}{|c|ccccc|}
      \hline
      $(K,\eps)$ & $N$ & $T_a$(sec) & $T_d$(sec) & $T_d / T_a$ & $\eps_a$\\
      \hline
      (2048,1e-4) & 7.82e+4 & 2.00e+1 & 3.87e+3 & 1.94e+2 & 1.15e-4\\
      (8192,1e-4) & 3.13e+5 & 8.80e+1 & 6.17e+4 & 7.02e+2 & 1.21e-4\\
      (32768,1e-4)& 1.25e+6 & 3.90e+2 & 9.90e+5 & 2.54e+3 & 1.07e-4\\
      \hline
      (2048,1e-6) & 7.82e+4 & 3.20e+1 & 3.87e+3 & 1.21e+2 & 1.04e-6\\
      (8192,1e-6) & 3.13e+5 & 1.38e+2 & 6.20e+4 & 4.50e+2 & 9.65e-7\\
      (32768,1e-6)& 1.25e+6 & 6.05e+2 & 1.01e+6 & 1.67e+3 & 1.20e-6\\
      \hline
      (2048,1e-8) & 7.82e+4 & 4.70e+1 & 3.87e+3 & 8.24e+1 & 8.58e-9\\
      (8192,1e-8) & 3.13e+5 & 2.03e+2 & 6.22e+4 & 3.06e+2 & 1.69e-8\\
      (32768,1e-8)& 1.25e+6 & 8.78e+2 & 9.95e+5 & 1.13e+3 & 1.33e-8\\
      \hline
    \end{tabular}
  \end{center}
  \caption{Results of an airfoil with the Helmholtz kernel.}
  \label{tbl:foil}
\end{table}

\begin{table}[h]
  \begin{center}
    \includegraphics[height=2in]{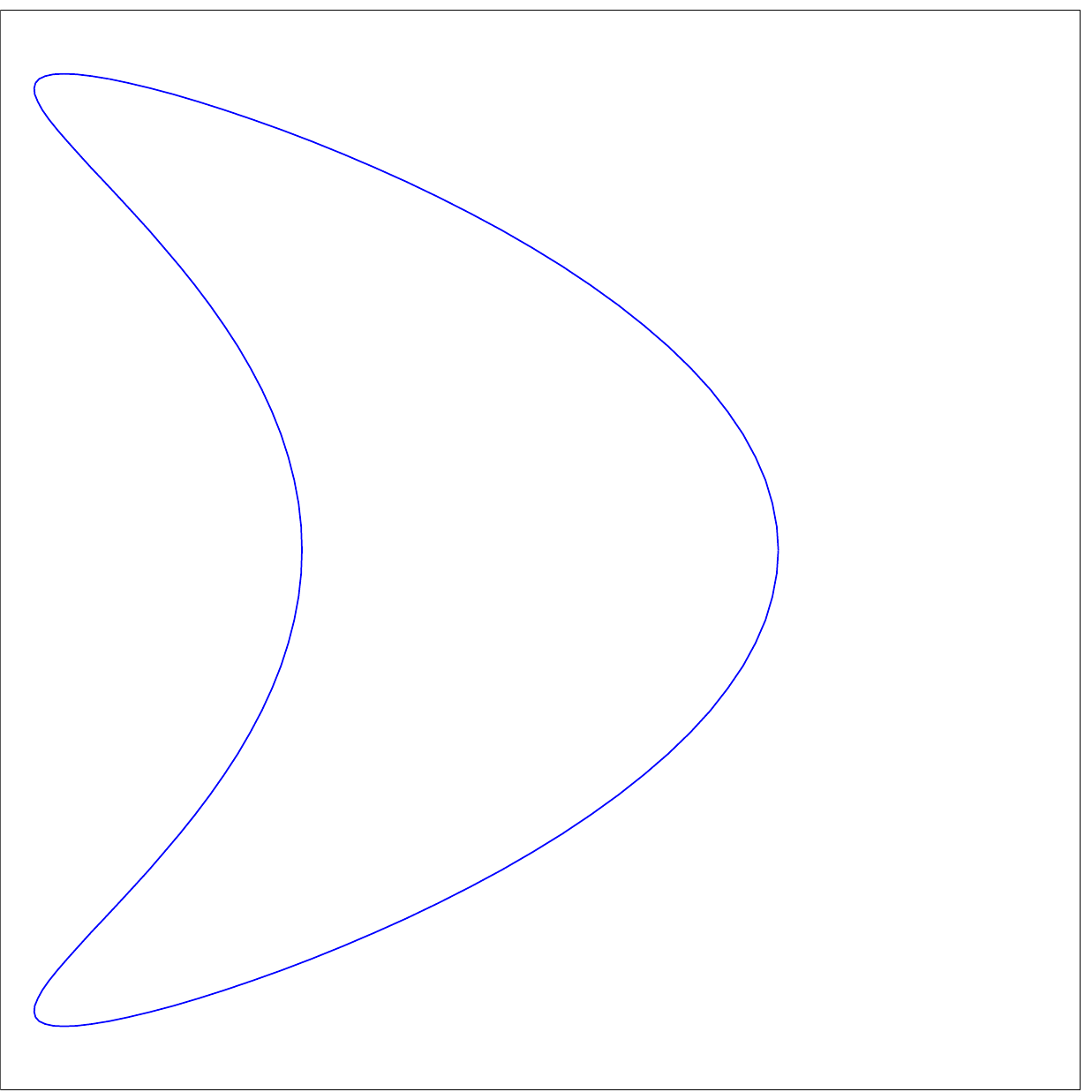}
    \includegraphics[height=2in]{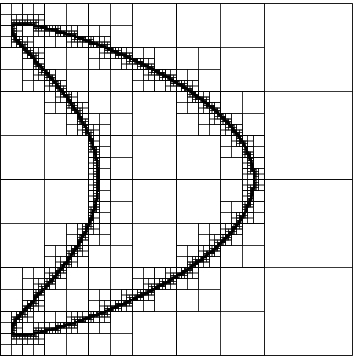}\\
    \vspace{0.1in}
    \begin{tabular}{|c|ccccc|}
      \hline
      $(K,\eps)$ & $N$ & $T_a$(sec) & $T_d$(sec) & $T_d / T_a$ & $\eps_a$\\
      \hline
      (2048,1e-4) & 1.13e+5 & 4.00e+1 & 8.11e+3 & 2.03e+2 & 1.08e-4\\
      (8192,1e-4) & 4.53e+5 & 1.77e+2 & 1.30e+5 & 7.36e+2 & 1.33e-4\\
      (32768,1e-4)& 1.81e+6 & 8.04e+2 & 2.09e+6 & 2.60e+3 & 1.41e-4\\
      \hline
      (2048,1e-6) & 1.13e+5 & 6.10e+1 & 8.11e+3 & 1.33e+2 & 9.35e-7\\
      (8192,1e-6) & 4.53e+5 & 2.72e+2 & 1.30e+5 & 4.78e+2 & 9.15e-7\\
      (32768,1e-6)& 1.81e+6 & 1.24e+3 & 2.10e+6 & 1.70e+3 & 8.80e-7\\
      \hline
      (2048,1e-8) & 1.13e+5 & 9.20e+1 & 8.16e+3 & 8.87e+1 & 1.45e-8\\
      (8192,1e-8) & 4.53e+5 & 4.05e+2 & 1.30e+5 & 3.22e+2 & 1.31e-8\\
      (32768,1e-8)& 1.81e+6 & 1.80e+3 & 2.11e+6 & 1.17e+3 & 1.52e-8\\
      \hline
    \end{tabular}
  \end{center}
  \caption{Results of a kite-shaped model with the Helmholtz kernel.}
  \label{tbl:kite}
\end{table}

Compared with the results presented in \cite{cheng-2006-riwfmm2d}, our
algorithm is slower by a factor of 8. The reason is that we heavily
use the kernel evaluation formula in our algorithm. The 2D Helmholtz
kernel involves the Hankel functions and the current computational
procedure for their evaluation is rather slow. On the other hand, all
of the high frequency translations in \cite{cheng-2006-riwfmm2d} are
precomputed and stored in the diagonal form and no special function
evaluation is required during the computation.

\begin{figure}
  \begin{center}
    \includegraphics[height=2in]{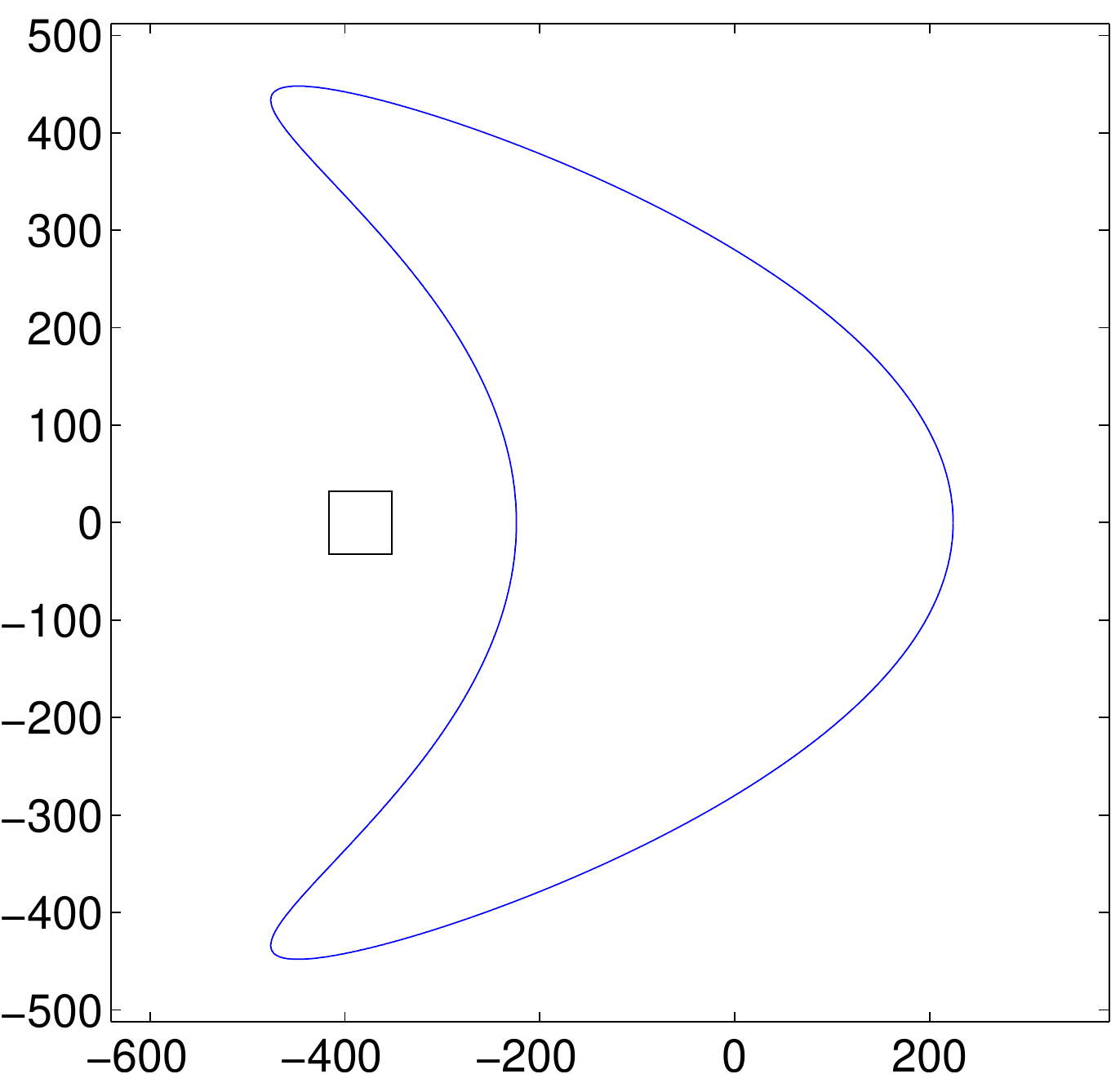} 
    \includegraphics[height=3in]{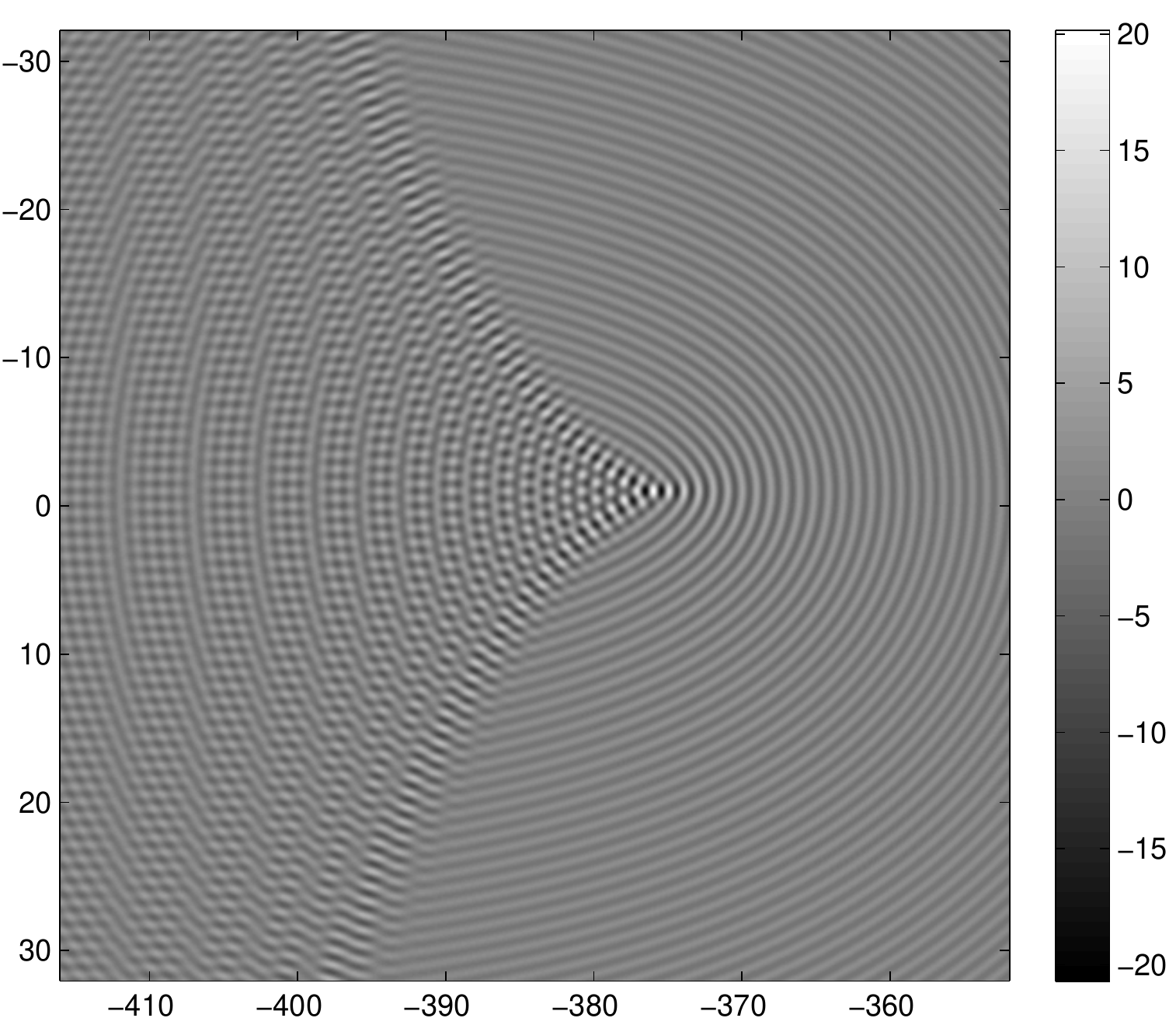}
  \end{center}
  \caption{Scattering field of the kite-shaped object with $K=1024$.
    Top: a square region that contains the caustics. Bottom: the real
    part of the scattering field inside the square. The field is
    sampled at 8 points per wavelength.  }
  \label{fig:kite}
\end{figure}

Finally, we apply our algorithm to the solution of the BIE formulation
\[
\frac{1}{2} \phi(x) + \int_{\p D} \left( \frac{\p G(x,y)}{\p n(y)} - \i \eta G(x,y) \right) \phi(y) d y
= - u^{inc}(x)
\]
of the 2D scattering problem mentioned in Section \ref{sec:intro}.
Here, we report the numerical results for the smooth objects in Tables
\ref{tbl:circ} and \ref{tbl:kite}. In our experiments, we use a
uniform discretization of about 20 points per wavelength. We pick
$\eta = \pi$ and set the incoming field $u^{inc}(x)$ to be $e^{2\pi\i
  x \cdot d}$ with $d=(1,0)$. We discretize the integral equation with
the Nystr\"{o}m method \cite{colton-1983-iemst,kress-1999-lie} and use
the endpoint-corrected trapezoidal rules from \cite{kapur-1997-hctqr}
to integrate the weakly singular part of the integral. The system is
solved iteratively using the GMRES algorithm and the restarted number
is set to be $80$. Within each iteration of the GMRES solver, the
application of the integral operator is accelerated using our
multidirectional algorithm with $\eps=$1e-4. Table \ref{tbl:sccirc}
summarizes the results for the circle with wavelengths from 1024 to
8192. Here $T_i$ is the averaged time of each iteration, $N_i$ is the
number of iterations, and $T_t$ is the total time. Table
\ref{tbl:sckite} reports the results of the kite-shaped object in
Table \ref{tbl:kite}.  In Figure \ref{fig:kite}, we display the
scattering field of the kite-shaped object in a region with caustics.

\begin{table}[h]
  \begin{center}
    \begin{tabular}{|c|cccc|}
      \hline
      $K$ & $N$ & $T_i$(sec) & $N_i$ & $T_t$(sec) \\
      \hline
      1024 & 65536 & 22 & 72 & 1.60e+3\\
      2048 & 131072 & 45 & 93 & 4.32e+3\\
      4096 & 262144 & 99 & 118 & 1.20e+4\\
      8192 & 524288 & 202 & 150 & 3.12e+4\\
      \hline
    \end{tabular}
  \end{center}
  \caption{
    Timings of computing the scattering field of the circle.
    $K$ is the size of the problem in terms of the wavelength,
    $N$ is the number of quadrature points,
    $T_i$ is the averaged time of each iteration, $N_i$ is the number of
    iterations, and $T_t$ is the total time.
  }
  \label{tbl:sccirc}
\end{table}

\begin{table}[h]
  \begin{center}
    \begin{tabular}{|c|cccc|}
      \hline
      $K$ & $N$ & $T_i$(sec) & $N_i$ & $T_t$(sec) \\
      \hline
      1024 & 65536 & 22 & 227 & 5.11e+3\\
      2048 & 131072 & 46 & 314 & 1.49e+4\\
      4096 & 262144 & 99 & 435 & 4.42e+4\\
      8192 & 524288 & 204 & 604 & 1.25e+5\\
      \hline
    \end{tabular}
  \end{center}
  \caption{
    Timings of computing the scattering field of the kite-shaped object.
  }
  \label{tbl:sckite}
\end{table}

\section{Conclusions}
\label{sec:concl}

In this paper, we described a directional multiscale algorithm for
computing the $N$-body problem for the high frequency Helmholtz kernel
in two dimensions. The approach follows the framework described in
\cite{engquist-2007-fdmaok}. Our algorithm is accurate and works well
for problems in all scales. By using the directional low rank
representations for regions that follow the directional parabolic
separation condition, our algorithm achieves the optimal $O(N\log N)$
complexity. A new and more efficient randomized technique compared to
the one in \cite{engquist-2007-fdmaok} has also been introduced for
the construction of the low rank separated representations. The
numerical results have shown that our algorithm is capable of
addressing very large scale problems in high frequency scattering.

For future work, we would like to have a rigorous proof for the
randomized procedure proposed in Section \ref{sec:direct}. Another
interesting direction for future research is to apply this kind of
directional multiscale idea to other problems with oscillatory
behavior, in both two and three dimensions.  One typical example is
the computation of the far field pattern of a scattering field
\cite{colton-1983-iemst, ying-2007-sftba}.

{\bf Acknowledgments.} The authors would like to thank P.G. Martinsson
for helpful discussions. B.E. is partially supported by an NSF grant
DMS 0714612 and a startup grant from the University of Texas at
Austin. L.Y. is partially supported by an Alfred P. Sloan Research
Fellowship and a startup grant from the University of Texas at Austin.

\bibliographystyle{abbrv}
\bibliography{ref}

\end{document}